\magnification=\magstep1
\baselineskip=15pt
\mathsurround=1pt

\font\Bbb=msbm10

\font\chapfont=cmssbx10 scaled \magstep3
\font\sectfont=cmssbx10 scaled \magstep2
\font\subsectfont=cmssbx10 scaled \magstep1
\font\exampfont=cmssbx10
\font\scrpt=cmsy10

\def\chap#1\par{\noindent{\chapfont#1}\par\medskip}
\def\sect#1\par{\noindent{\sectfont#1}\par\nobreak\medskip}
\def\subsect#1\par{\bigbreak\noindent{\subsectfont#1}\par\nobreak\medskip}
\def\pf#1:{\noindent{\bf#1}:}
\def\examp#1. {\medbreak\noindent{\exampfont#1}. }

\def\emph{\it}

\def\arrow-under #1 {\buildrel #1 \over \longrightarrow} 
\def\longarrow-under #1 {\buildrel #1 \over \loongrightarrow}
\def\longerarrow-under #1 {\buildrel #1 \over \looongrightarrow}
\def\larrow-under #1 {\buildrel #1 \over \longleftarrow}

\def\Sec˜{\mathhexbox278}
\def\Ñ{\kern-1pt\vrule height2.4ptwidth3ptdepth-1.9pt\kern1pt}
\def\¾{\hfill \rlap{$\sqcup$}$\sqcap$\par\bigskip}

\def\Zü{\hbox{\Bbb Z}}
\def\QÎ{\hbox{\Bbb Q}}
\def\Rå{\hbox{\Bbb R}}
\def\C'{\hbox{\Bbb C}}

\def\scrptè{{\cal U}}

\def\p{\pi}

\def\s{\sigma}

\def\c{\raise 1pt\hbox{$\chi$}}

\def\G{\Gamma}

\def\bdyÆ{\partial}
\def\elÉ{\ell}
\def\infÔ{\infty}
\def\empÑ{\hbox{\Bbb ?}}

\def\simil°{\sim}
\def\homot¦{\simeq}
\def\iso¤{\approx}

\def\º{\equiv}
\def\­{\neq}
\def\Î{\in}
\def\Ï{\notin}
\def\É{\cdots}
\def\²{\leq}
\def\³{\geq}
\def\:{\,\colon}
\def\—{\pm}

\def\unionÚ{\cup}
\def\UnionÚ{\bigcup}
\def\intÛ{\cap}
\def\IntÛ{\bigcap}
\def\DisjÜ{\coprod}
\def\disjÜ{\amalg}
\def\cross´{\times}

\def\lbr›{\{}
\def\rbrÍ{\}}
\def\vbar|{\mid}
\def\res|{\rlap{\raise1pt\hbox{$|$}}\lower1pt\hbox{$|$}}
\def\<ø{\langle}
\def\>ù{\rangle}

\def\to¬{\rightarrow}
\def\To¬{\longrightarrow}
\def\mpstoØ{\mapsto}
\def\incl«{\hookrightarrow}
\def\subsetofÞ{\subset}

\def\barŠ  {\overline}
\def\hatš  {\widehat}
\def\tild'  {\widetilde}

\def\Diff{\hbox{\it Diff\/}}
\def\rel{{\rm \, rel\ }}
\def\scrK{{\cal K}}
\def\scrM{{\cal M}}
\def\scrE{{\cal E}}

\centerline{\sectfont Spaces of Knots}
\medskip 
\centerline{Allen Hatcher}
\bigskip

Classical knot theory is concerned with isotopy classes of knots in the 3\Ñsphere, in other words, path-components of the space $ \scrK $ of all smooth submanifolds of $ S^3 $ diffeomorphic to the circle $ S^1 $. What can be said about the homotopy types of these various path-components? One would like to find, for the path-component $ \scrK_K $ containing a given knot $ K $, a small subspace $ \scrM_K $ to which $ \scrK_K $ deformation retracts, thus a minimal homotopic model for $ \scrK_K $. In this paper we describe a reasonable candidate for $ \scrM_K $ and prove for many knots $ K $ that $ \scrK_K $ does indeed have the homotopy type of the model $ \scrM_K $. The proof would apply for all $ K $ provided that a certain well-known conjecture in 3\Ñmanifold theory is true, the conjecture that every free action of a finite cyclic group on $ S^3 $ is equivalent to a standard linear action.

The model $ \scrM_K $ takes a particularly simple form if $ K $ is either a torus knot or a hyperbolic knot. In these cases $ \scrM_K $ is a single orbit of the action of $ SO(4) $ on $ \scrK_K $ by rotations of the ambient space $ S^3 $, namely an orbit of a ``maximally symmetric" position for $ K $, a position where the subgroup $ G_K \subsetofÞ SO(4) $ leaving $ K $ setwise invariant is as large as possible. The orbit is thus the coset space $ SO(4)/G_K $. The assertion that $ \scrK_K $ deformation retracts to $ \scrM_K = SO(4)/G_K $ is then a sort of {\emph homotopic rigidity\/} property of $ K $.

Let us describe these models in more detail. Consider first the case of the trivial knot. Its most symmetric position is clearly a great circle in $ S^3 $. The subgroup $ G_K \subsetofÞ SO(4) $ taking $ K $ to itself is then the index two subgroup of $ O(2) \cross´ O(2) $ consisting of orientation-preserving isometries. It was shown in [H1] that $ \scrK_K $ has the homotopy type of the orbit $ \scrM_K = SO(4)/G_K $, which can be identified with the $ 4 $\Ñdimensional Grassmann manifold of $2$\Ñplanes through the origin in $ \Rå^4 $.

Consider next a nontrivial torus knot $ K = K_{p,q} $ for relatively prime integers $ p $ and $ q $, neither of which is $ \—1 $. Regarding $ S^3 $ as the unit sphere in $ \C'^2 $, the most symmetric position for $ K $ is as the set of points $ (z^p, z^q)/\sqrt 2 $ with $ |z| = 1 $. Its symmetry group then contains the unitary diagonal matrices with $ z^p $ and $ z^q $ as diagonal entries, forming a subgroup $ S^1 \subsetofÞ SO(4) $. There is also a rotational symmetry reversing the orientation of $ K $, given by complex conjugation in each variable. Thus $ G_K $ contains a copy of $ O(2) $, whose restriction to $ K $ is the usual action of $ O(2) $ on $ S^1 $. It is easy to see that $ G_K $ cannot be larger than this since if it were, $ K $ would be pointwise fixed by a nontrivial element of $ SO(4) $, hence would be unknotted. We will show that $ \scrK_K $ has the homotopy type of the orbit $ SO(4)/G_K $, a closed $ 5 $\Ñmanifold. This may have been known for some time, but there does not seem to be a proof in the literature.

Now let $ K $ be hyperbolic, so $ S^3 - K $ has a unique complete hyperbolic structure, and let $ \G_K $ be the finite group of orientation-preserving isometries of this hyperbolic structure. An easy argument in hyperbolic geometry shows that elements of $ \G_K $ must take meridians of $ K $ to meridians, so the action of $ \G_K $ on $ S^3 - K $ extends to an action on $ S^3 $. By the Smith conjecture [MB], no nontrivial elements of $ \G_K $ fix $ K $ pointwise, so $ \G_K $ is a group of diffeomorphisms of $ K $, hence $ \G_K $ must be cyclic or dihedral. Assuming the action of $ \G_K $ on $ S^3 $ is equivalent to an action by elements of $ SO(4) $, then we can isotope $ K $ to a ``symmetric" position in which the action is by isometries of $ S^3 $, so we have an embedding $ \G_K \subsetofÞ SO(4) $. We will show in this case that $ \scrK_K $ has the homotopy type of $ SO(4)/\G_K $. The symmetry group $ G_K $ cannot be larger than $ \G_K $, so in this symmetric position for $ K $ we have $ G_K = \G_K $, and $ SO(4)/\G_K $ is the orbit of $ K $ under the $ SO(4) $ action.

The hypothesis that the action of $ \G_K $ on $ S^3 $ is equivalent to an isometric action can be restated as the well-known conjecture that the orbifold $ S^3/\G_K $ has a spherical structure, a special case of Thurston's geometrization conjecture for orbifolds. We call it {\emph the linearization conjecture for\/} $ K $. It is a theorem of Thurston that the geometrization conjecture is true for orbifolds which are not actually manifolds. So the linearization conjecture for $ K $ is true unless $ \G_K $ is a cyclic group acting freely on $ S^3 $. It is also known to be true for free actions by a cyclic group of order two [L], a power of two [M], or three [R2]. It appears to be still unknown for prime orders $ p \³ 5 $.

To prove these results we study the space $ \scrE $ of smooth embeddings $ f \: S^1 \to¬ S^3 $. The space $ \scrK $ is the orbit space of $ \scrE $ under the action of the diffeomorphism group $ \Diff(S^1) $ by composition in the domain. This is a free action, and the projection $ \scrE \to¬ \scrK $ is a principal bundle with fiber $ \Diff(S^1) $. The path-components of $ \scrE $ are the oriented knot types, since a standard orientation of $ S^1 $ induces an orientation on the image of each embedding $ S^1 \to¬ S^3 $. A component $ \scrK_K $ of $ \scrK $ is the image of one or two components of $ \scrE $, depending on whether $ K $ is invertible or not.

The group $ SO(4) $ also acts on $ \scrE $, by composition in the range. This is a free action when $ K $ is nontrivial, defining a principal bundle $ \scrE_K \to¬ \scrE_K/SO(4) $ for each component $ \scrE_K $ of $ \scrE $. We show that $ \scrE_K/SO(4) $ is a $ K(\p,1) $ for all nontrivial $ K $, and we give a description of the group $ \p = \p_K $. In particular we can deduce that $ \scrE_K/SO(4) $ has the homotopy type of a finite CW complex, and it follows that the same is true for $ \scrE_K $. By contrast, we have only been able to show that $ \scrK_K $ has finite homotopy type if we assume the linearization conjecture.

When $ K $ is a nontrivial torus knot the group $ \p_K $ is trivial, so $ \scrE_K/SO(4) $ is contractible and $ \scrE_K $  has the homotopy type of a single orbit $ SO(4) $. When $ K $ is a hyperbolic knot, $ \p_K $ is $ \Zü $, and it follows that $ \scrE_K/SO(4) \homot¦ S^1 $ and $ \scrE_K \homot¦ S^1 \cross´ SO(4) $. If the linearization conjecture is true for $ K $, $ \scrE_K $ has the homotopy type of a single orbit of the action of $ SO(2) \cross´ SO(4) $ on $ \scrE_K $ by rotations in both domain and range, namely the orbit containing a sufficiently symmetric embedding. Such an orbit has the form $ (SO(2) \cross´ SO(4))/\G^+_K $ where $ \G^+_K $ is the cyclic subgroup of $ \G_K $ consisting of symmetries preserving the orientation of $ K $. 

Knots which are not torus knots or hyperbolic knots are satellite knots, and for these the situation becomes more complicated. In particular the homotopic rigidity property of torus knots and hyperbolic knots fails for satellite knots. Modulo the linearization conjecture again, we show that $ \scrK_K $ has the homotopy type of a model $ \scrM_K $ which is a finite-dimensional manifold of the form $ (X_K \cross´ SO(4))/\G_K $ where $ \G_K $ is a finite group of ``supersymmetries" of $ K $ and $ X_K $ is the product of a torus of some dimension and a number of configuration spaces $ C_n $ of ordered $ n $-tuples of distinct points in $ \Rå^2 $. These configuration spaces occur only when the satellite structure of $ K $ involves nonprime knots. When they are present, $ \p_1 \scrK_K $ involves braid groups, a phenomenon observed first in [G] in the case that $ K $ itself is nonprime. When there are no configuration spaces, $ \scrM_K $ is a closed manifold, but in general there can be no closed manifold model for $ \scrK_K $ since it may not satisfy Poincar\'e duality. The space $ X_K $ appearing in $ \scrM_K $ is determined just by the general form of the satellite structure of $ K $, while the group $ \G_K $ is more delicate, depending strongly on the particular knots appearing in the satellite structure. ($ \G_K $ is the quotient of $ \p_0\Diff^+(S^3 \rel K) $ by the subgroup generated by Dehn twists along essential tori in $ S^3 - K $.) See Section 2 of the paper, where the case of satellite knots is treated in detail.

\bigskip

{\bf\noindent 1. Homotopically Rigid Knots}
\medskip

As described in the introduction, let $ \scrE_K $ and $ \scrK_K $ be the components of the spaces of embeddings $ S^1 \to¬ S^3 $ and images of such embeddings, respectively, corresponding to a given knot $ K $. These two spaces are related via the fibration that defines $ \scrK_K $ as a quotient space of $ \scrE_K $,

$$
F \To¬ \scrE_K \To¬ \scrK_K
$$
whose fiber $ F $ is the diffeomorphism group $ \Diff(S^1) $ if $ K $ is invertible, or the orientation-preserving subgroup $ \Diff^+(S^1) $ if $ K $ is not invertible. Here ``invertible" has its standard meaning of ``isotopic to itself with reversed orientation." Since $ \Diff(S^1) \homot¦ O(2) $ and $ \Diff^+(S^1) \homot¦ SO(2) $, the homotopy types of $ \scrE_K $ and $ \scrK_K $ should be closely related. It happens that $ \scrE_K $ is more directly accessible to our techniques, so we study this first, then apply the results to $ \scrK_K $.

When $ K $ is nontrivial, the group $ SO(4) $ acts freely on $ \scrE_K $ by composition in the range, defining a principle bundle
$$
SO(4) \To¬ \scrE_K \To¬ \scrE_K/SO(4)
$$

\proclaim Theorem 1. If $ K $ is nontrivial, $ \scrE_K/SO(4) $ is aspherical, i.e., a $ K(\p,1) $. Its fundamental group $ \p $ is trivial if $ K $ is a torus knot, and $ \Zü $ if $ K $ is hyperbolic. Hence $ \scrE_K \homot¦ SO(4) $ if $ K $ is a torus knot, and $ \scrE_K \homot¦ S^1 \cross´ SO(4) $ if $ K $ is hyperbolic.

\pf Proof: By restriction of orientation-preserving diffeomorphisms of $ S^3 $ to a chosen copy of the knot $ K $ we obtain a fiber bundle
$$
\Diff^+(S^3\rel K) \To¬ \Diff^+(S^3) \To¬ \scrE_K
$$
where ``rel $K $" indicates diffeomorphisms which restrict to the identity on $ K $. The fiber bundle property is a special case of the general result that restriction of diffeomorphisms to a submanifold defines a fiber bundle; see [L]. When we factor out the action of $ SO(4) $ on $ \Diff^+(S^3) $ and $ \scrE_K $ by composition in the range we obtain another fiber bundle
$$
\Diff^+(S^3\rel K) \To¬ \Diff^+(S^3)/SO(4) \To¬ \scrE_K/SO(4)
$$
By the Smale Conjecture [H1], the total space $ \Diff^+(S^3)/SO(4) $ of this bundle has trivial homotopy groups, hence $ \p_i(\scrE_K/SO(4)) \iso¤ \p_{i+1}\Diff^+(S^3\rel K) $ for all $ i $. In fact, since the bundle is a principal bundle with contractible total space, $ \scrE_K/SO(4) $ is a classifying space for the group $ \Diff^+(S^3\rel K) $.

In similar fashion, if $ N $ is a tubular neighborhood of $ K $ in $ S^3 $ we have a fibration
$$
\Diff(S^3\rel N) \To¬ \Diff^+(S^3\rel K) \To¬ E(N\rel K)
$$
where $ E(N\rel K) $ is the space of embeddings $ N \incl« S^3 $ restricting to the identity on $ K $. It is a standard fact that $ E(N\rel K) $ has the homotopy type of the space of automorphisms of the normal bundle of $ K $ in $ S^3 $. Since $ K $ is diffeomorphic to $ S^1 $, $ E(N\rel K) $ thus has the homotopy type of $ \Zü \cross´ S^1 $. The $ \Zü $ factor measures different choices of a nonzero section of the normal bundle. Since elements of $ \Diff^+(S^3\rel K) $ must take longitude to longitude, up to isotopy, for homological reasons, we may as well replace the base space $ E(N\rel K) $ of the bundle by $ S^1 $.

Let $ M $ be the closure of $ S^3 - N $, a compact manifold with torus boundary. We can then identify $ \Diff(S^3\rel N) $ with $ \Diff(M\rel\bdyÆ M) $. It has been known since the early 1980's that  $ \Diff(M\rel\bdyÆ M) $ has contractible components, as a special case of more general results about Haken manifolds [H2],[I]. Therefore to show that $ \Diff^+(S^3\rel K) $ has contractible components, hence that $ \scrE_K/SO(4) $ is aspherical, it suffices to verify that in the long exact sequence of homotopy groups for preceding the fibration, the boundary map from  $ \p_1E(N\rel K) $ to $ \p_0\Diff(M\rel\bdyÆ M) $ is injective. 

A generator for the infinite cyclic group $ \p_1E(N \rel K) $ is represented by a full rotation of each disk fiber. Under the boundary map this gives a diffeomorphism of $ M $ supported in a collar neighborhood of $ \bdyÆ M $ which restricts to a standard Dehn twist in each meridional annulus. One can see this is nontrivial in $ \p_0\Diff(M\rel\bdyÆ M) $ by looking at the induced homomorphism on $ \p_1(M,x_0) $ for a basepoint $ x_0 \Î \bdyÆ M $. Namely, it is conjugation by a meridian loop, so if the Dehn twist were isotopic to the identity fixing $ \bdyÆ M $, this conjugation would be the trivial automorphism of $ \p_1(M,x_0) $ and hence the meridian would lie in the center of $ \p_1(M,x_0) $. However, the only Haken manifolds with nontrivial center are the orientable Seifert fiberings, with the fiber generating the center [BZ]. But the only knot complements which are Seifert-fibered are torus knots, with the fiber being non-meridional. Thus the Dehn twist is nontrivial in $ \p_0\Diff(M\rel\bdyÆ M) $, and the same argument shows that any nonzero power of it is also nontrivial, so the boundary map is injective.

The next thing to show is that $ \p_0\Diff^+(S^3 \rel K) = 0 $ if $ K $ is a torus knot, or equivalently that $ \p_0\Diff(M\rel\bdyÆ M) $ is generated by meridional Dehn twists. This can be shown by standard 3\Ñmanifold techniques, in the following way. The manifold $ M $ is Seifert fibered over a disk with two multiple fibers of distinct multiplicities. In particular, $ M $ is the union of two solid tori intersecting along an annulus $ A $. Modulo meridional Dehn twists, a diffeomorphism of $ M $ fixing $ \bdyÆ M $ can be isotoped, staying fixed on $ \bdyÆ M $, so that it takes $ A $ to itself. The restriction of the diffeomorphism to $ A $ must be isotopic to the identity $ \rel \bdyÆ A $ since it extends over the solid tori. The argument is completed by appealing to the fact that $ \p_0\Diff(T \rel \bdyÆ T) = 0 $ for $ T $ a solid torus.

Suppose now that $ K $ is hyperbolic. From the long exact sequence of homotopy groups for the restriction fibration $ \Diff(M) \to¬ \Diff(\bdyÆ M) $ we obtain a short exact sequence
$$
0 \To¬ \p_1\Diff(\bdyÆ M) \arrow-under {\bdyÆ} \p_0\Diff(M \rel \bdyÆ M) \To¬ \p_0\Diff_0(M) \To¬ 0
$$
where $ \Diff_0(M) $ consists of the diffeomorphisms of $ M $ whose restriction to $ \bdyÆ M $ is isotopic to the identity. Injectivity of the map $ \bdyÆ $ was shown two paragraphs above, since the image of this map is generated by longitudinal and meridional Dehn twists near $ \bdyÆ M $. By famous theorems of Waldhausen and Cerf, along with Mostow rigidity, we have $ \p_0\Diff(M) \iso¤ Isom(M) $, the finite group of hyperbolic isometries of $ M $. The group $ \p_0\Diff_0(M) $ is a subgroup of this, the isometries whose restriction to the cusp torus is a rotation. Since hyperbolic isometries are locally determined, each isometry in $ \p_0\Diff_0(M) $ is uniquely determined by its restriction rotation of $ \bdyÆ M $.

We can see that $ \p_0\Diff(M \rel \bdyÆ M) $ is isomorphic to $ \Zü \cross´ \Zü $ by the following argument from [HM]. The fiber $ \Diff(M \rel \bdyÆ M) $ of the map $ \Diff(M) \to¬ \Diff(\bdyÆ M) $ has the same homotopy groups as the homotopy fiber, whose points are pairs consisting of a diffeomorphism of $ M $ together with an isotopy of its restriction to $ \bdyÆ M $ to the identity. Thus we can view $ \p_0\Diff(M \rel \bdyÆ M) $ as rotations of $ \bdyÆ M $ with a homotopy class of isotopies of the rotation to the identity. This is the same as a lift of the rotation to a translation of the universal cover $ \Rå^2 $. Thus we can identify $ \p_0\Diff(M \rel \bdyÆ M) $ with a group of translations of $ \Rå^2 $ containing the group $ \p_1\Diff(\bdyÆ M) $ of deck transformations as a finite-index subgroup. Hence $ \p_0\Diff(M \rel \bdyÆ M) \iso¤ \Zü \cross´ \Zü $.

To show that $ \p_0\Diff^+(S^3 \rel K) \iso¤ \Zü $ we need to see that a meridional Dehn twist of $ M $ generates a $ \Zü $ direct summand of $ \p_0\Diff(M \rel \bdyÆ M) $. First note that no hyperbolic isometry of $ M $ can restrict to a purely meridional rotation of $ \bdyÆ M $ since such an isometry would extend to a periodic diffeomorphism of $ S^3 $ fixing $ K $ pointwise, which is ruled out by the Smith conjecture [MB]. Thus a meridional Dehn twist of $ M $ is not a proper multiple of any other element of $ \p_0\Diff(M \rel \bdyÆ M) \iso¤ \Zü \cross´ \Zü $. The meridional Dehn twist is then one element of a basis for $ \p_0\Diff(M \rel \bdyÆ M) $, which is to say that it generates a direct summand.

This finishes the proof of the statements about $ \scrE_K/SO(4) $. We deduce the statements about $ \scrE_K $ by looking at the bundle $ \scrE_K \to¬ \scrE_K/SO(4) $. If $ K $ is a torus knot we have shown that the base space of this bundle has trivial homotopy groups, hence is contractible, so the total space deformation retracts onto a fiber. If $ K $ is hyperbolic the base space has the homotopy type of $ S^1 $, so the total space has the homotopy type of a principal $ SO(4) $ bundle over $ S^1 $. Since $ SO(4) $ is connected, this bundle must be trivial. \¾

When $ K $ is hyperbolic, a generator for the $ \Zü $ summand of $ \p_1\scrE_K $ can be described as follows. A generator of the cyclic group $ \p_0\Diff_0(M) $ is represented by a diffeomorphism which extends to a periodic diffeomorphism $ g \: (S^3,K) \to¬ (S^3,K) $ preserving orientations of both $ S^3 $ and $ K $. Let $ g_t \: S^3 \to¬ S^3 $ be an isotopy from the identity to $ g = g_1 $. Restricting $ g_t $ to $ K $ gives a path of embeddings $ f_t \: K \to¬ S^3 $ with $ f_1(K) = f_0(K) = K $. The restriction of $ f_1 $ to $ K $ is a diffeomorphism $ K \to¬ K $ isotopic to the identity, say by an isotopy $ h_t $, so if we follow the isotopy $ f_t $ by the inverse isotopy $ h_{1-t} $, we obtain a loop of embeddings of $ K $ in $ S^3 $. This generates the $ \Zü $ summand of $ \p_1\scrE_K $. 

If $ K $ satisfies the linearization conjecture, as described in the introduction, then the diffeomorphism $ g $ is conjugate, via an element $ \s \Î \Diff^+(S^3) $, to an element of $ SO(4) $. Replacing $ K $ by the equivalent knot $ \s(K) $, we may assume that $ g $ itself is in $ SO(4) $. For the isotopy $ g_t $ we can then choose an arc in an $ S^1 $ subgroup of $ SO(4) $ containing $ g $. The isotopy $ f_t $ in the preceding paragraph is then constant in $ \scrE_K/SO(4) $, so a loop generating $ \p_1(\scrE_K/SO(4)) $ consists of the reparametrizations $ h_t $ of $ K $. This implies that $ \scrE_K $ has the homotopy type of the coset space $ (S^1 \cross´ SO(4))/\G^+_K $ where $ \G^+_K $ is the cyclic group generated by $ (g\res| K,g^{-1}) $. 

Factoring out parametrizations, which have the homotopy type of $ O(2) $ or $ SO(2) $ depending on whether $ K $ is invertible or not, we obtain:

\proclaim Corollary. If $ K $ is a torus knot, or a hyperbolic knot satisfying the linearization conjecture, then $ \scrK_K $ has the homotopy type of $ SO(4)/\G_K $. \¾

For a hyperbolic knot $ K \subsetofÞ S^3 $ the symmetry group $ G_K $ is always a subgroup of $ \G_K $. This is because the natural map $ G_K \to¬ \p_0\Diff^+(S^3-K) \iso¤ \G_K $ is injective, by the theorem that a periodic diffeomorphism of a Haken manifold which is homotopic to the identity must be part of an $ S^1 $ action (see [T],[FY]), which happens only when $ S^3 - K $ is a Seifert manifold, hence $ K $ is a torus knot.

\bigskip
{\bf\noindent 2. Satellite Knots }
\medskip
(This section is not yet written.)

\bigskip\noindent
{\bf References}
\medskip\parindent=33pt

\item{[BZ]} G. Burde and H. Zieschang, Eine Kennzeichnung der Torusknoten, Math. Ann. 167 (1966), 169-176.

\item{[FY]} M. Freedman and S.-T. Yau, Homotopically trivial symmetries of Haken manifolds are toral, Topology 22 (1983), 179-189.

\item{[G]} A. Gramain, Sur le groupe fondamental de l'espace des noeuds, Ann. Inst. Fourier 27, 3 (1977), 29-44.

\item{[H1]} A. Hatcher, A proof of the Smale conjecture,  Ann. of Math. 117 (1983), 553-607.

\item{[H2]} A. Hatcher,  Homeomorphisms of sufficiently large $ P^2$\Ñirreducible
3\Ñmanifolds.\break Topology 15 (1976), 343-347. For a more recent version, see the paper
``Spaces of incompressible surfaces" available on the author's webpage: 

http://math.cornell.edu/$\,\tilde{}\,$hatcher

\item{[HM]} A. Hatcher and D. McCullough, Finiteness of classifying spaces of relative diffeomorphism groups of 3\Ñmanifolds, Geometry and Topology 1 (1997), 91-109.

\item{[I]} N. V. Ivanov,  Diffeomorphism groups of
Waldhausen manifolds, J. Soviet Math. 12 (1979), 115-118  (Russian
original in Zap. Nauk. Sem. Leningrad Otdel. Mat.  Inst. Steklov 66
(1976) 172-176.  Detailed write-up:  Spaces of surfaces in
Waldhausen manifolds, Preprint LOMI P-5-80 Leningrad (1980)).

\item{[L]} E. L. Lima, On the local triviality of the restriction map for embeddings, Comment. Math. Helv. 38 (1964), 163-164.

\item{[L]} G. R. Livesay, Fixed point free involutions on the 3\Ñsphere, Ann. of Math. 72 (1960), 603-611.

\item{[MB]} J. W. Morgan and H. Bass, eds., {\it The Smith Conjecture}, Academic Press, 1984.

\item{[M]} R. Myers, Free involutions on lens spaces, Topology 20 (1981), 313-318.

\item{[R1]} J. H. Rubinstein, Free actions of some finite groups on $ S^3 $, Math. Ann. 240 (1979), 165-175.

\item{[R2]} J. H. Rubinstein, unpublished.

\item{[T]} J. Tollefson, Homotopically trivial periodic homeomorphisms of 3\Ñmanifolds, Ann. of Math. 97 (1973), 14-26.

\end